\input BRIAN.STY
\def\[{{[\![}}
\def\]{{]\!]}}
\emergencystretch3mm

\MSC{20g35,20f28,20f40,17b20,17b40}
\UDC{512.544.43+512.743+512.554.37}

\title{Automorphisms and isomorphisms of Chevalley groups and algebras}
\author{Anton A. Klyachko}
\address{\myAddress}
\grants{\RFBR 05-01-00895}
\leftline{\it Dedicated to Rebecca}\bigskip

\abstract{\narrower
An adjoint Chevalley group of rank at least 2
over a rational algebra (or a similar ring), its
elementary subgroup,\break and the corresponding Lie
ring have the same automorphism group. These automorphisms are explicitly
described.
}

\s
0. Introduction

Suppose that $\Phi$ is a reduced irreducible root system,
$R$ is an associative commutative ring with unity, $G(\Phi,R)$ is the
corresponding adjoint Chevalley group, and $E(\Phi,R)$ is its elementary
subgroup (see Section 5).

There are a lot of results (see, e.g., [Wat80], [Pet82], [GMi83],
[HO'M89], [Abe93], [Che00], [Bun07], and references therein%
\fn{%
Unfortunately, some interesting papers on this subject
(e.g., [Abe93])
contain errors.
}) asserting
that, under some conditions, all automorphisms of Chevalley (or similar)
groups are standard in some sense (depending on what a particular author
succeeded to prove). In this paper, we use the most universal and natural
definition of standardness suggested by A.~E.~Zalesskii [Zal83]: an
automorphism of an adjoint Chevalley group is called \emph{standard} if it
is induced by an automorphism of the corresponding Lie algebra. More
precisely, this means the following. Clearly, $E(\Phi,R)$ and $G(\Phi,R)$
embed naturally into the automorphism group of the corresponding Lie
algebra~$L(\Phi,R)$ over $R$. A slightly less obvious fact is that (under
some conditions, see Section~5) both groups are normal in
$\Aut_R L(\Phi,R)$ and even in the larger group
$\Aut_\Z L(\Phi,R)=\Aut_\Z R\semitimes \Aut_R L(\Phi,R)$ consisting of the
automorphisms of this algebra considered as a Lie ring.  Thus, each
automorphism $f\in\Aut_\Z L(\Phi,R)$ of the Lie ring induces an
automorphism $f'\:g\mapsto fgf^{-1}$ of the Chevalley groups $G(\Phi,R)$
and $E(\Phi,R)$.  The main results of this paper are the following
theorems.

\proclaim{Automorphism theorem}.
For any reduced irreducible root system $\Phi$ of rank $\ge2$, there
exists an integer $m$ such that, for any associative commutative ring $R$
without additive torsion, with unity and $1\over m$, all automorphisms of
the Chevalley group $G(\Phi,R)$ and its elementary subgroup
$E(\Phi,R)$ are standard; the groups
$\Aut_\Z L(\Phi,R)$, $\Aut_\Z R\semitimes \Aut_R L(\Phi,R)$,
$\Aut G(\Phi,R)$, and $\Aut E(\Phi,R)$ are
isomorphic; the map
$\Aut_\Z L(\Phi,R)\ni f\mapsto f'\in \Aut G(\Phi,R)$ is a group
isomorphism; a similar map $\Aut_\Z L(\Phi,R)\to\Aut E(\Phi,R)$ is a group
isomorphism also.

\proclaim{Isomorphism theorem}.
For any reduced irreducible root system $\Phi$ of rank $\ge2$, there
exists an integer $m$ such that, for any associative commutative rings
$R$ and $R'$
without additive torsion,
with unity and $1\over m$,
there are natural one-to-one correspondences between the
following three sets:
$$
\matrix{
\{\hbox{group isomorphisms } G(\Phi,R)\to G(\Phi,R')\},
\quad
\{\hbox{group isomorphisms } E(\Phi,R)\to E(\Phi,R')\},
\cr
\hbox{and}\quad
\{\hbox{Lie ring isomorphisms } L(\Phi,R)\to L(\Phi,R')\},
\cr
}
$$
i.e., each group isomorphism $G(\Phi,R)\to G(\Phi,R')$ maps $E(\Phi,R)$
onto $E(\Phi,R')$; each group isomorphism $E(\Phi,R)\to E(\Phi,R')$ can be
extended uniquely to an isomorphism $G(\Phi,R)\to G(\Phi,R')$;
each ring isomorphism $f\:L(\Phi,R)\to L(\Phi,R')$ induces a group
isomorphism
$\Aut_\Z L(\Phi,R)\supseteq G(\Phi,R)\to
G(\Phi,R')\subseteq \Aut_\Z L(\Phi,R')$ by the formula
$\phi\mapsto f\phi f^{-1}$; each group isomorphism
$G(\Phi,R)\to G(\Phi,R')$ is induced by such a way by a unique ring
isomorphism.
\newline
\hbox{\quad }Each ring isomorphism $f\:L(\Phi,R)\to L(\Phi,R')$ is
semilinear, i.e., $f(rx)=\alpha(r)f(x)$ for some ring isomorphism
$\alpha\: R\to R'$ uniquely determined by $f$.

In particular, these theorems allow us to describe the
automorphisms of all adjoint Chevalley groups of rank at least 2
over any commutative $\Q$-algebras.
Similar results were obtained by Yu Chen [Che95], [Che96] (see also
[Che00]), but he assumed additionally that $R$ is an algebra over $\Q$
without zero divisors.

The idea of describing the automorphisms of linear groups by the passage to
the Lie algebras was first introduced and applied by
V.~M.~Levchuk [Lev83] and
E.~I.~Zelmanov
[Zel85]. We use this general idea, but our approach is quite different.

The above theorems reduce the problem of finding
automorphisms/isomorphisms of Chevalley groups to an (easier) analogous
problem for Chevalley algebras. The automorphisms of Chevalley algebra
are explicitly described in Section~7. Each automorphism of
$L(\Phi,R)$ is a composition of an inner automorphism (i.e., a conjugation
by an element of $G(\Phi,R)$) and automorphisms induced by symmetries of
the corresponding Dynkin diagram.

Our proofs are completely calculation-free and use only few properties of
Chevalley groups. Thus, this approach can work in a more general setting.
An \emph{elementary group scheme} $E$ is a subgroup of
$\SL_n(\Z[z_1,z_2,\dots])$ generated by some matrices
$\{x_i(z_j)\;;\; i\in I,\;j=1,2,\dots\}$.  For an elementary group scheme
$E$, the symbol $E(R)$ denotes the subgroup of $\SL_n(R)$ consisting of
all matrices of the form $a(r_1,r_2,\dots)$, where $a\in E$ and
$r_j\in R$.  We say that $E(R)$ is an $n$-dimensional \emph{$R$-group}.
Clearly, $E(R)$ is generated by the matrices
$\{x_i(r)\;;\; i\in I,\; r\in R\}$.
For an $R$-group $E(R)$ we consider the following conditions:
\item{(EX)}
\emph{Exponentiality}:
$x_i(z_1)x_i(z_2)=x_i(z_1+z_2)$ for all $i\in I$.

\item{(AL)} \emph{Algebraicity}:
$E(R[t])$ is a normal subgroup of a linear algebraic group
$G\subseteq\SL_n(R[t])$ defined by some polynomial equations with integer
coefficients. The group $E(R[t])$ is the normal closure of its subgroup
$E(R)$.

\item{(PC$_S$)} \emph{Power conjugacy}:
two matrices $x_i$ and $x_i^s$ are conjugate in $E(R)$ for each
$i\in I$ and each $s\in S$, where $S\subseteq\Z$ is a set of integers.

\Example.
In Section 5, we show that, under the assumption of the above theorems,
an adjoint elementary Chevalley group
$E(\Phi,R)$
has properties (EX), (AL),
and (PC$_S$), where $S=\Z\cap\{a^2\;;\;a\in R^*\}$.

The author is very grateful to N. A. Vavilov for valuable remarks.
The author also thanks 
V.~M.~Levchuk,
A.~V.~Mikhalev, 
V.~A.~Petrov, 
and 
D.~A.~Timashev
for reading this paper and useful discussion.

\s
1. Nullstellensatz

Recall that an ideal is called \emph{radical} if the corresponding factor
ring has no nonzero nilpotents. We use the following form of Hilbert's
Nullstellensatz.

\proclaim{Nullstellensatz}.
Suppose that $g,f_1,\dots,f_l\in\Z[y_1,\dots,y_m]$ are some polynomials
and the quasi-identity
$$
\forall r_1,\dots,r_l\in R \quad
f_1(r_1,\dots)=0\ \&\ \dots\ \&\ f_l(r_1,\dots)=0\ \imp\
g(r_1,\dots)=0
$$
holds for $R=\C$.
Then there exists positive integer $b$ such that the quasi-identity
$$
\forall r_1,\dots,r_l\in R \quad
f_1(r_1,\dots)=0\ \&\ \dots\ \&\ f_l(r_1,\dots)=0\ \imp\
(g(r_1,\dots))^b=0
$$
holds for any associative commutative ring $R$ with unity and without
additive torsion. If the ideal of $\Z[y_1,\dots,y_m]$ generated by
$f_1,\dots,f_l$ is radical, then we can take $b=1$.

\s
2. Unipotence

A matrix $A$ is called \emph{unipotent} if $A-1$ is a nilpotent matrix.
We say that an automorphism of an $R$-group $E(R)$ is \emph{unipotent} if
it sends all $x_i(r)$ to unipotent matrices.
An automorphism $\phi$ is said to be \emph{$m$-unipotent} if
$(\phi(x_i(r))-1)^m=0$ for all $i\in I$ and $r\in R$.

\Proposition 1.
Suppose that an associative commutative ring $R$ with unity has no
additive torsion. Then, for any integers $n\ge 1$, $p\ge2$, and $d\ge1$,
there exist positive integers $q$ and $m$ such that any automorphism of an
$n$-dimensional $R$-group with Property {\rm(PC$_{\{p,q^d\}}$)} is
$m$-unipotent.

\Proof
If the $R$-group satisfies Property {\rm(PC$_{\{p,q^d\}}$)}, then, for any
automorphism $\phi$, the matrices $\phi(x_i(r))$, $\phi(x_i(r))^p$, and
$\phi(x_i(r))^{q^d}$ are conjugate. Thus, Proposition 1 is a corollary of
the following lemma.

\Lemma 1.
For any integers $n\ge 1$, $p\ge2$, and $d\ge1$, there exist positive
integers $q$ and $m$ such that, if the characteristic polynomials of
matrices $A, A^p, A^{q^d}\in\SL_n(R)$ over an associative commutative ring
$R$ with unity and without additive torsion coincide, then $(A-1)^m=0$.

\Proof
First, assume that $R=\C$. Then, these three matrices have the same set of
eigenvalues and
raising to the power $p$ acts as a permutation of these eigenvalues.
Hence, for any eigenvalue $\lambda$,

$$
\lambda^{p^{n!}-1}=1
\quad\hbox{and, by the same reason,}\quad
\lambda^{q^{dn!}-1}=1.
$$
Clearly, these equalities imply $\lambda=1$ if we take, e.g., $q=p^{n!}-1$
(then $p^{n!}-1$ and $q^{dn!}-1$ are coprime).
So, the assertion is proven for the case $R=\C$.

The condition
$$
\hbox{
$\det A=1$ and
the characteristic polynomials of $A$, $A^p$, and $A^{q^d}$ coincide}
\eqno{(*)}
$$
is a system of integer-coefficient polynomial equations on
entries of matrix $A$.
Each complex root $B\in M_n(\C)$ of this system is unipotent
(if $q$ is chosen as above).
By Nullstellensatz,
this implies that, if a matrix $A$ over $R$
satisfies $(*)$, then each entry $c_{ij}$ of the matrix $C=(A-1)^n$
satisfies the equality $c_{ij}^b=0$ for some integer $b$. Hence,
$C^{bn^2}=(A-1)^{bn^3}=0$. This completes the proof of Lemma 1 and
Proposition 1.

\s
3. Curves

Take an $R$-group $E(R)$. We say that the group $E(R[t])$ is the
\emph{group of curves on the group $E(R)$}.

Clearly, for any curve $g(t)\in E(R[t])$ and any polynomial $f(t)\in R[t]$,
the curve $\REP_f(g)\:=g(f(t))$ also belongs to $E(R[t])$. We say that
$g(f(t))$ is the \emph{reparametrisation of the curve $g$ by means of the
polynomial $f$}. Thus, $\REP_f\:E(R[t])\to E(R[t])$ is an endomorphism of
the group of curves.

Performing fairly standard calculations
$$
(1+tX+t^2Y+o(t^2))^{-1}=
1-(tX+t^2Y)+(tX+t^2Y)^2+o(t^2)=
1-tX+t^2(X^2-Y)+o(t^2),
$$
$$
\eqalign{
&[(1+tX_1+t^2Y_1+o(t^2)),(1+tX_2+t^2Y_2+o(t^2))]=\cr
&=(1+tX_1+t^2Y_1+o(t^2))(1+tX_2+t^2Y_2+o(t^2))
(1+tX_1+t^2Y_1+o(t^2))^{-1}(1+tX_2+t^2Y_2+o(t^2))^{-1}=\cr
&=1+t^2(Y_1+Y_2-Y_1-Y_2+X_1^2+X_2^2+X_1X_2-X_1^2-X_1X_2-X_2X_1-X_2^2+X_1X_2)
+o(t^2)=\cr
&=1+t^2(X_1X_2-X_2X_1)+o(t^2),
}
$$
we obtain a fairly standard formula:
$$
[(1+tX_1+o(t)),(1+tX_2+o(t))]=\REP_{t^2}(1+t\[X_1,X_2\]+o(t))+o(t^2),
\eqno{(1)}
$$
where $[x,y]\:=xyx^{-1}y^{-1}$ is the group commutator and
$\[x,y\]\:=xy-yx$ is the ring commutator.

\s
4. Continuity and smoothness

The set of matrices
$
T(E(R))=\big\{X\in M_n(R)\;|\; 1+tX+t^{2}Y\in E(R[t])
\hbox{ \ for some $Y\in M_n(R[t])$}\big\}
$
\smallskip\noindent
is called the \emph{tangent module} of an $R$-group $E(R)$.
Clearly, this set is an $R[E(R)]$-module, i.e., it is closed with respect
to
\-addition: $(1+tX+o(t))(1+tY+o(t))=1+t(X+Y)+o(t)$;
\-multiplication by scalars:
$\REP_{rt}(1+tX+o(t))=1+tXr+o(t)$;
\-the action of the group $E(R)$:
  $g(1+tX+o(t))g^{-1}=(1+tgXg^{-1}+o(t))$ (In what follows, we put
  $g\o X\:=gXg^{-1}$).

\noindent
If the tangent module is a Lie algebra, i.e., if it is closed with respect
to ring commutator $\[A,B\]=AB-BA$, we call this module the \emph{tangent
algebra}. We say that an $n$-dimensional $R$-group $E(R)$ is
\emph{adjoint} if $T(E(R))$ is a Lie algebra isomorphic as an
$R[E(R)]$-module to $R^n$ (with the natural action of $E(R)$).

We say that an automorphism $\phi$ of an $R$-group $E(R)$ is
\emph{quasicontinuous} if it can be extended to an automorphism $\~\phi$
of the group of curves such that $\~\phi$ commutes with all
integer-coefficient reparametrisations:
$\~\phi(\REP_f(g))=\REP_f(\~\phi(g))$ for all $g\in E(R[t])$ and all
$f\in\Z[t]$. The automorphism $\phi$ is called \emph{continuous} if it is
quasicontinuous, the automorphism~$\~\phi$ is quasicontinuous,
the automorphism $\kern1pt\~{\kern-2pt\~\phi}$ is quasicontinuous, and so
on (infinitely many times).

Put $E_k(R)\:=E(R[t])\cap(1+t^kM_n(R[t]))$. Since $\ker\REP_0=E_1(R)$, we
have the equality $\~\phi(E_1(R))=E_1(R)$ for any continuous automorphism
$\phi$.  We say that a continuous automorphism $\phi$ is \emph{smooth}
(two times differentiable) if $\~\phi(E_k(R))=E_k(R)$ for $k=1,2,3$. Note
that the continuity [smoothness] of an automorphism implies the
continuity [smoothness] of the inverse automorphism. In Section 5, we
show that any continuous automorphism of a Chevalley group is smooth
(under some conditions).

\Proposition 2.
Any smooth automorphism $\phi$ of a group $E(R)$ induces an automorphism
$d\phi$ (the \emph{differential of $\phi$}) of the tangent
module considered as an abelian group. In addition, we have
$$
\~\phi(1+tX+o(t))=1+td\phi(X)+o(t)
\quad\hbox{and}\quad
d\phi(g\o X)=\phi(g)\o d\phi(X)
\quad\hbox{for all $g\in E(R)$ and $X\in T(E(R))$}.
$$
If $T(E(R))$ is a Lie algebra, then $d\phi$ is an automorphism of this
algebra considered as a Lie ring.

\Proof
If $X\in T(E(R))$, then $1+tX+t^2Y\in E(R[t])$ for some
$Y\in M_n(R[t])$ and $\~\phi(1+tX+t^{2}Y)=1+tZ+o(t)$
for some $Z\in M_n(R)$ (because $E_1(R)$ is an invariant subgroup).
Put $d\phi(X)=Z$.
This is well defined, because $\~\phi$ leaves invariant the subgroup
$E_2(R)$. The bijectivity of $d\phi$ follows from the smoothness of
$\phi^{-1}$.
The equalities
$$
\matrix{
\~\phi((1+tX+o(t))(1+tY+o(t)))&=&
\~\phi(1+t(X+Y)+o(t))=1+td\phi(X+Y)+o(t)\hfill\cr
\cr
||&\cr
\cr
\~\phi(1+tX+o(t))\~\phi(1+tY+o(t))&=&
(1+td\phi(X)+o(t))(1+td\phi(Y)+o(t))=
1+t(d\phi(X)+d\phi(Y))+o(t)\cr
}
$$
show that $d\phi$ is an endomorphism of the additive group.
A similar argument
$$
\matrix{
\~\phi(g(1+tX+o(t))g^{-1})&=&
\~\phi(1+tgXg^{-1}+o(t))=1+td\phi(gXg^{-1})+o(t)\hfill\cr
\cr
||&\cr
\cr
\phi(g)\~\phi(1+tX+o(t))\phi(g)^{-1}&=&
\phi(g)(1+td\phi(X)+o(t))\phi(g)^{-1}=
1+t\phi(g)d\phi(X)\phi(g)^{-1}+o(t)\cr
}
$$
proves the equality
$d\phi(g\o X)=\phi(g)\o d\phi(X)$.

The automorphism $\~\phi$ commutes with integer-coefficient
reparametrisations, leaves invariant $E_3(R)$, and, hence, maps equality
(1) to
$$
[\~\phi(1+tX_1+o(t)),\~\phi(1+tX_2+o(t))]=
\REP_{t^2}\~\phi(1+t\[X_1,X_2\]+o(t))+o(t^2).
$$
Therefore,
$$
[(1+td\phi(X_1)+o(t)),(1+td\phi(X_2)+o(t))]=
\REP_{t^2}(1+td\phi(\[X_1,X_2\])+o(t))+o(t^2).
$$
Applying formula (1) to the left-hand side, we obtain
$$
\eqalign{
\REP_{t^2}(1+t\[d\phi(X_1),d\phi(X_2)\]+o(t))+o(t^2)&=
[(1+td\phi(X_1)+o(t)),(1+td\phi(X_2)+o(t))]=\cr
&=\REP_{t^2}(1+td\phi(\[X_1,X_2\])+o(t))+o(t^2).
}
$$
Thus, $\[d\phi(X_1),d\phi(X_2)\]=d\phi(\[X_1,X_2\])$.
This proves that $d\phi$ is an endomorphism of the tangent algebra.

\medskip

If $E(R)$ is adjoint, then it embeds naturally into the automorphism group
$\Aut_{\Z}T(E(R))$ of its tangent algebra considered as a Lie ring.

\Proposition 3.
Any smooth automorphism of an adjoint $R$-group $E(R)$ is standard, i.e.,
it has the form $\phi(g)=\alpha g\alpha^{-1}$, where $\alpha$ is an
automorphism of the Lie ring $T(E(R))$ normalising the subgroup $E(R)$.

\Proof
This follows immediately from Proposition 2, we can take $\alpha=d\phi$.

\Proposition 4.
Suppose that a commutative associative ring $R$ with unity and
${1\over q!}$ has no additive torsion, an $R$-group $E(R)$ has properties
{\rm(EX)} and {\rm(AL)}, and $\phi$ and $\phi^{-1}$ are mutually inverse
$q$-unipotent automorphisms of $E(R)$. Then these automorphisms are
continuous.

\Proof
Take a matrix $a(t)\in E(R[t])$. Clearly, $a(r)\in E(R)$ for any $r\in R$.
Let us prove that
$$
\hbox{\sl the matrix $\phi(a(k))$ depends polynomially on the number
$k\in\Z$,}
$$
i.e., there exists a matrix $b_{a}(t)\in \SL_n(R[t])$ such
that $\phi(a(k))=b_{a}(k)$ for all $k\in\Z$. (Note that the absence of
additive torsion implies the uniqueness of such matrix $b_{a}(t)$.)

Indeed, it is sufficient to prove this for $a(t)=x_i(rt^l)$, because these
matrices generate the group $E(R[t])$. So,
$$
\phi(x_i(rk^l))=(\phi(x_i(r)))^{k^l}
\quad\hbox{by Property {\rm(EX)}}.
$$
But $(\phi(x_i(r)))^m$ depends polynomially on $m$, as the matrix
$\phi(x_i(r))$ is unipotent:
$$
(\phi(x_i(r)))^m=(1+A)^m=1+mA+{m(m-1)\over2}A^2+\dots+
{m(m-1)\dots(m-q+1)\over q!}A^q.
$$
Thus, we can extend the automorphism $\phi$ to the group $E(R[t])$
putting $\~\phi(a(t))\:=b_{a}(t)$.

Let us prove that $\~\phi(a(t))=b_{a}(t)$ lies in $E(R[t])$.
For each integer $k$, the matrix $b_{a}(k)$ belongs to $E(R)$ and, hence,
belongs to the group $G$ defined by integer-coefficient polynomial
equations (see Property (AL)). Therefore, the matrix~$b_{a}(t)$
satisfies the same equations. Thus,
$b_{a}(t)\in G$ and we have
$$
\matrix{
&\hbox to 5.8mm{$E(R)$=\hss}
&\~\phi(E(R))&        \subseteq  &E(R[t])\cr
\cr
&&\lower-1pt\hbox{$\scriptstyle|$}\cap&                 &\lower-1pt\hbox{$\scriptstyle|$}\cap\cr
\cr
&&\~\phi(E(R[t]))&\subseteq &G\cr
}
\quad\hbox{and}\quad
\matrix{
E(R[t])=\nc{E(R)}_{E(R[t])}\subseteq&\nc{E(R)}_G&\supseteq
\nc{E(R)}_{\tilde{\phi}(E(R[t]))}=\~\phi(E(R[t]))\cr
\cr
&\parallel&\hbox{(by Property (AL))}\cr
\cr
&E(R[t]),&\cr
}
$$
where $\nc X_H$ means the normal closure of a set $X$ in a group $H$.
Thus, $\~\phi(E(R[t]))\subseteq E(R[t])$.

The automorphism $\phi^{-1}$ also can be extended to the group of curves
and $\~\phi\~{(\phi^{-1})}(a(k))=\~{(\phi^{-1})}\~\phi(a(k))=a(k)$ for any
$k\in\Z$ and any $a(t)\in E(R[t])$. Clearly, this implies the equalities
$\~\phi\~{(\phi^{-1})}(a(t))=\~{(\phi^{-1})}\~\phi(a(t))=a(t)$
(because $R$ has no additive torsion) and the bijectivity of $\~\phi$.

By the construction, the automorphism $\~\phi$ commutes
with all integer-coefficient reparametrisations. So, $\phi$ is
quasicontinuous.
Clearly, $\~\phi$ is also $q$-unipotent and, hence, quasicontinuous.
Thus, an obvious induction argument completes the proof of the
continuity of $\phi$.

\s
5. Chevalley groups

Suppose that $\Phi$ is a reduced irreducible root system, $L(\Phi)$ is the
corresponding simple complex Lie algebra. The algebra~$L(\Phi)$ has a
basis $h_1, h_2,\dots, x_1, x_2,\dots$ (the \emph{Chevalley basis}) such
that the structure constants are integer and the matrices of the operators
$(\ad x_i)^k/k!$ are integer and nilpotent for all $k\in\N$. The
\emph{Chevalley algebra} is the Lie $R$-algebra~$L(\Phi,R)$ with the same
structure constants.

Suppose that $N(\Phi)=\Aut_\C L(\Phi)$ is the automorphism group of the
algebra $L(\Phi)$, and $G(\Phi)=(\Aut_\C L(\Phi))^\o$ is the connected
component of the identity of this group. The algebraic groups
$G(\Phi)\subseteq N(\Phi)\subseteq\GL(L(\Phi))\subset\SL_n(\C)$ are
defined over $\Z$. Let $R$ be an associative commutative ring with unity
and let $N(\Phi,R)$ and $G(\Phi,R)$ be the groups of $R$-rational points
of $N(\Phi)$ and $G(\Phi)$, i.e., the subgroups of $\SL_n(R)$ (where
$n=1+\dim L(\Phi)$) defined by the same integer-coefficient polynomial
equations as the groups $N(\Phi)$ and $G(\Phi)$, respectively, (in the
Chevalley basis). Note that $N(\Phi,R)=\Aut_R L(\Phi,R)$, because the
property of being an automorphism can be written as a system of
integer-coefficient polynomial equations (depending on the structure
constants). The group $G(\Phi,R)$ is called the \emph{\(adjoint\)
Chevalley group}.  The group $E(\Phi,R)\subseteq G(\Phi,R)$ generated by
the matrices $x_i(r)=\exp(\ad rx_i)$, where $r\in R$, is called the
\emph{elementary subgroup} of the Chevalley group~$G(\Phi,R)$.

\Example.
For the root system $A_l$, we have $L(A_l)=\Sl_{l+1}(\C)$ is the Lie
algebra consisting of all traceless matrices, $L(A_l,R)=\Sl_{l+1}(R)$,
$G(A_l,R)=\PGL_{l+1}(R)$, and $E(A_l,R)=\PE_{l+1}(R)$ is the subgroup of
$\PGL_{l+1}(R)$ generated by the images of the transvections $1+rE_{ij}$,
where $i\ne j$ and $r\in R$. (Note that, for some rings, this group
$\PGL_{l+1}(R)$ can be large than the central quotient of the general
linear group $\GL_{l+1}(R)$.)

In the following lemma, we summarise some (probably) known properties of
Chevalley groups and algebras.

\Lemma 2.
Let $\Phi$ be a reduced irreducible root system of rank $\ge2$ and let
$R$ be an associative commutative ring without additive torsion,
with unity and $1\over6$. Then
\item{\rm(i)}
the group $E(\Phi,R)$ is an $R$-group with Properties {\rm(EX)} and
{\rm(PC$_S$)}, where $S=\Z\cap\{a^2\;;\;a\in R^*\}$;
\item{\rm(ii)}
for each subgroup $H$ of $G(\Phi,R)$ normalised by $E(\Phi,R)$, there
exists a unique ideal $J$ of $R$ such that $H$ is contained in
$G(\Phi,R)\cap(1+M_n(J))$ and contains the normal closure
$\nc{\{x_i(r)\;;\; r\in J\}}_{E(\Phi,R)}$ of the set
$\{x_i(r)\;;\; r\in J\}$;
\item{\rm(iii)}
$E(\Phi,R)$ is an automorphism invariant (i.e., characteristic)
subgroup of $G(\Phi,R)$;
\item{\rm(iv)}
$E(\Phi,R)$ satisfies Property {\rm(AL)};
\item{\rm(v)}
$\Aut_\Z L(\Phi,R)\iso\Aut_\Z R\semitimes \Aut_R L(\Phi,R)$;
\item{\rm(vi)}
in the group $\Aut_\Z L(\Phi,R)$, the subgroups $G(\Phi,R)$ and
$E(\Phi,R)$ are normal and their centralisers are trivial.

\Proof

\item{\rm(i)}{%
Property (EX) follows immediately from the definition.  Steinberg's
relation R5, $h_i(s)x_i(r)h_i(s)^{-1}=x_i(s^2r)$ (see, e.g., [VPl96]),
where $r\in R$, $s\in R^*$, and $h_i(s)\in E(\Phi,R)$ are some particular
matrices, shows that Property~(PC$_S$) holds too.
}

\item{\rm(ii)}{%
Taking into account that $G(\Phi,R)$ is centreless in the adjoint case
[AHu88], we see that (ii) is a slightly weakened form of the well-known
theorem on subgroups of Chevalley groups normalised by the
elementary subgroups [Vas86] (see also [ASu76], [Abe89], [Gol97], [CKe99],
[VGN06]).
}

\item{\rm(iii)}{%
This was also proven by Vaserstein in [Vas86]. Note that, in [HaV03],
it was in fact proven the endomorphism invariance of the elementary
subgroup of a Chevalley group.
}

\item{\rm(iv)}{%
The normality of $E(\Phi,R[t])$ in the linear algebraic group
$G(\Phi,R[t])$ defined by polynomial equations with integer coefficients
follows immediately from (iii).  The equality
$\nc{E(\Phi,R)}_{E(\Phi,R[t])}=E(\Phi,R[t])$ follows from (ii). Indeed,
put $H=\nc{E(\Phi,R)}_{E(\Phi,R[t])}$. The inclusion
$E(\Phi,R)\subseteq H\subseteq G(\Phi,R[t])\cap(1+M_n(J))$ implies
$J=R[t]$. Therefore,
$E(\Phi,R[t])=\nc{\{x_i(f)\;;\; f\in R[t]\}}_{E(\Phi,R[t])}=
\nc{\{x_i(f)\;;\; f\in J\}}_{E(\Phi,R[t])} \subseteq H$ and
$H=E(\Phi,R[t])$.
}

\item{\rm(v)}{%
Let $U$ be the algebra $L(\Phi,R)$ considered as a left module over
itself. Then
$$
\End_{L(\Phi,R)}U=R,
\quad\hbox{i.e., all endomorphisms are scalar multiples of the identity.}
\eqno{(**)}
$$
Indeed, this is true for $R=\C$, because the algebra $L(\Phi,\C)$ is
simple. Therefore, $(**)$ holds for any $R$ without additive torsion,
because both conditions on a matrix, being an endomorphism of $U$ and
being a scalar multiple of the identity, are integer-coefficient
systems of linear equations on the entries of the matrix.
\newline
Note that $(**)$ remains valid if we consider $U$ as a module over Lie
ring $L(\Phi,R)$, i.e., each endomorphism $f$ of $U$ must be $R$-linear.
Indeed, for any $u\in U$, there exist $y_i\in L(\Phi,R)$ and $u_i\in U$
such that $u=\sum(\ad y_i)(u_i)$, because
$L(\Phi,R)=[L(\Phi,R),L(\Phi,R)]$. Therefore,
$$
ru=\sum(\ad ry_i)(u_i)
\quad\hbox{and}\quad
f(ru)=f\(\sum(\ad ry_i)(u_i)\)=\sum(\ad ry_i)f(u_i)=
r\sum(\ad y_i)f(u_i)=rf(u).
$$
Now,
take an automorphism $\phi$ of the ring $L(\Phi,R)$ and
consider the algebra $L(\Phi,R)$ as an $L(\Phi,R)$-module $U_\phi$
with action $(y,u)\mapsto (\ad\phi(y))u$. The mapping
$u\mapsto\phi(u)$ is an isomorphism between
the modules $U$ and $U_\phi$ over the Lie ring $L(\Phi,R)$. This
isomorphism induces an isomorphism of endomorphism rings
$R=\End_{L(\Phi,R)}U
\mathop{\longrightarrow}\limits^{\alpha_\phi}\End_{L(\Phi,R)}U_\phi=R$.
Thus, we have a homomorphism $\Aut_\Z L(\Phi,R)\to\Aut_\Z R$,
$\phi\mapsto\alpha_\phi$, whose
kernel is $\Aut_R L(\Phi,R)$.
The right inverse homomorphism $\Aut_\Z R\to \Aut_\Z L(\Phi,R)$ maps
$\alpha\in\Aut_\Z R$ to the obvious automorphism of the Lie ring
$L(\Phi,R)=L(\Phi,\Z)\otimes R$ induced by $\alpha$.
So, we obtain the required decomposition of $\Aut_\Z L(\Phi,R)$ into the
semidirect product.
}

\item{\rm(vi)}{%
{\bf Normality.}
By virtue of (iii), it is sufficient to prove the normality of $G(\Phi,R)$.
For $R=\C$, this property is well known, see, e.g., [VOn88].
Let $F_N(y_{ij})=0$ and $F_G(y_{ij})=0$ be systems of integer-coefficient
polynomial equations that define the groups $N(\Phi,R)=\Aut_R L(\Phi,R)$
and $G(\Phi,R)$ (these systems do not depend on $R$).
We assume that the ideals of $\Z[y_{11},y_{12},\dots, y_{nn}]$ generated
by the sets of polynomials $F_N(y_{ij})$ and $F_G(y_{ij})$ are radical.
For $R=\C$ we have the quasi-identity
$$
F_G(Y)=0\ \&\ F_N(Z)=0\ \imp\ F_G(ZYZ^{-1})=0.
\eqno{(2)}
$$
Since the ideal of
$\Z[y_{11},y_{12},\dots, y_{nn},z_{11},z_{12},\dots, z_{nn}]$
generated by $F_G(Y)$ and $F_N(Z)$ is radical,
Nullstellensatz implies
that
quasi-identity (2) holds for all rings $R$ without additive
torsion. Thus, $G(\Phi,R)$ is a normal subgroup of $N(\Phi,R)$.
\newline
{\bf Centralisers.}
For $R=\C$, the centraliser of the set $\{x_i(1)\}$ in $\Aut_R L(\Phi,R)$
is trivial. Therefore, the same is true for any ring $R$ without additive
torsion (by Nullstellensatz). Thus, the centraliser of
the set $\{x_i(1)\}$ in the group
$\Aut_\Z L(\Phi,R)=(\Aut_\Z R)\semitimes\Aut_R L(\Phi,R)$
coincide with $\Aut_\Z R$.
On the other hand, each nontrivial ring automorphism $\alpha\in\Aut_\Z R$
induces a nontrivial automorphism $x_i(r)\mapsto x_i(\alpha(r))$ of
$E(\Phi,R)$. Therefore, the centraliser of $E(\Phi,R)$ in the group
$\Aut_\Z L(\Phi,R)$ is trivial. This completes the proof of Lemma 2.
}

\Proposition 5.
Let $\Phi$ be a reduced irreducible root system of rank $\ge2$ and let
$R$ be an associative commutative ring without additive torsion,
with unity and $1\over6$. Then any retraction
$\pi\:E(\Phi,R[t])\to E(\Phi,R)$ (i.e., a homomorphism such that $\pi^2=\pi$)
has the form $E(\Phi,R[t])\ni a(t)\mapsto a(r)\in E(\Phi,R)$ for some
$r\in R$. In other words, $\pi=\REP_r$.

\Proof
According to Lemma 2 (ii),
$$
\nc{\{x_i(f)\;;\; f\in J\}}_{E(\Phi,R[t])}\subseteq
\ker\pi\subseteq E(\Phi,R[t])\cap(1+M_n(J))
\quad\hbox{for some ideal $J$ of $R[t]$}.
$$
The right-hand inclusion and the equality
$E(\Phi,R[t])=E(\Phi,R)\semitimes\ker\pi$ show that $t-r\in J$ for some
$r\in R$; the left-hand inclusion and the equality
$E(\Phi,R)\cap\ker\pi=\1$ show that $J=(t-r)R[t]$.  Therefore,
$\ker\pi=E(\Phi,R[t])\cap(1+M_n(J))$ and $\pi=\REP_r$.

\medskip

Thus, we have a natural one-to-one correspondence between the ring $R$ and
the set of retractions. Clearly, the ring structure on $R$ can also be
described in terms of retractions and integer-coefficient
reparametrisations:
$$
\eqalign{
&\REP_{r+r'}\:E(\Phi,R[t])\To^{t\to t+t'}E(\Phi,R[t,t'])
\To^{t\to r}_{t'\to r'}E(\Phi,R),
\cr
&\REP_{rr'}\:E(\Phi,R[t])\To^{t\to tt'}E(\Phi,R[t,t'])
\To^{t\to r}_{t'\to r'}E(\Phi,R).
}
\eqno{(3)}
$$
Proposition 5 and these formulae imply that any continuous automorphism
$\phi\in\Aut E(\Phi,R)$ induces a ring automorphism $\^\phi\in\Aut_\Z R$
by the formula $\phi\REP_r\~\phi^{-1}=\REP_{\hat\phi(r)}$:
$$
\matrix{
E(\Phi,R[t]) & \To^{\~\phi} & E(\Phi,R[t])\cr
\downarrow\scriptstyle{\REP_r}&&\downarrow\scriptstyle{\REP_{\hat\phi(r)}}
\cr E(\Phi,R) & \To^{\phi}& E(\Phi,R).}
$$
For each ideal $J\nin R$ we have two normal subgroups of
$E(\Phi,R)$, namely, $E(J)_{\max}\:=E(\Phi,R)\cap(1+M_n(J))$ and
$E(J)_{\min}\:=\nc{\{x_i(r)\;;\; r\in J\}}_{E(\Phi,R)}$.

\Lemma 3.
Let $\Phi$ be a reduced irreducible root system of rank $\ge2$ and let $R$
be an associative commutative ring without additive torsion, with unity
and $1\over6$. Then
$\phi(E(J)_{\min})=E(\^\phi(J))_{\min}$ and
$\phi(E(J)_{\max})=E(\^\phi(J))_{\max}$
for any continuous automorphism
$\phi$ of the group $E(\Phi,R)$.

\Proof
Clearly, it is sufficient to define $E(J)_{\min}$ and $E(J)_{\max}$ in
terms of retractions. The subgroup
$E_1(\Phi,R)\:=E(\Phi,R[t])\cap(1+tM_n(R[t]))$ can be
defined as $E_1(\Phi,R)=\ker\REP_0$ (hence, this
subgroup is $\~\phi$-invariant). Then,
$$
E(J)_{\min}=
\nc{\{\REP_r(a(t))\;;\; r\in J,\ a(t)\in E_1(\Phi,R)\}}_{E(\Phi,R)}.
$$
The inclusion $\supseteq$ follows from the equality
$E_1(R)=
\nc{\{x_i(rt^k)\;;\;i\in I,\; r\in R,\; k=1,2\dots\}}_{E(\Phi,R[t])}$,
which is valid for any $R$-group with Property (EX).
$$
E(J)_{\max}=
\vtop{
\hsize 10cm
\noindent
the (unique) maximal subgroups among all
normal subgroups $H$ such that $E(J)_{\min}\subseteq H$ and
$E(J')_{\min}\not\subseteq H$ for any ideal $J'\not\subseteq J$.}
$$
The correctness of this definition of $E(J)_{\max}$ follows from Lemma 2
(ii) and the equality
$E(J_1+J_2)_{\min}=E(J_1)_{\min}\cdot\nobreak E(J_2)_{\min}$.

\Lemma 4.
Let $\Phi$ be a reduced irreducible root system of rank $\ge2$ and let $R$
be an associative commutative ring without additive torsion, with unity
and $1\over6$. Then
any  continuous automorphism $\phi$ of $E(\Phi,R)$ is smooth.

\Proof
We have to prove that the subgroups
$E_k(\Phi,R)\:=E(\Phi,R[t])\cap(1+t^kM_n(R[t]))$ are $\~\phi$-invariant.
This is true for $k=1$, because $E_1(\Phi,R)=\ker\REP_0$.
On the other hand, $E_1=E(tR[t])_{\max}$. Hence, the ideal
$tR[t]\nin R[t]$ is $\^{\~\phi\,}$-invariant by Lemma 3.
Then, the ideal $(tR[t])^k$ is $\^{\~\phi\,}$-invariant and the
subgroup $E_k(\Phi,R)=E((tR[t])^k)_{\max}$ is $\~\phi$-invariant.

\Proposition 6.
The tangent module of a Chevalley group coincide with the corresponding Lie
algebra: $T(E(\Phi,R))=L(\Phi,R)$.

\Proof
Suppose that $X\in T(E(\Phi,R))$, i.e., $1+tX+o(t)\in E(\Phi,R[t])$.
Let us express this element via the generators:
$$
1+tX+o(t)=\prod_j x_{i_j}(r_jt^{k_j})
\eqno{(4)}
$$
Clearly, we can assume that $k_j\in\{0,1\}$. Also, the substitution
$t=0$ shows that
$$
{\prod_j}' x_{i_j}(r_j)=1.
\quad\hbox{where the
prime means that the
product is taken over all $j$ such that $k_j=0$}.
$$
Therefore, the expression (4) can be rewritten in the form
$$
1+tX+o(t)=\prod_l g_lx_{i_l}(r_lt)g_l^{-1},
\quad\hbox{where $g_l\in E(\Phi,R)$.}
$$
Hence, $X=\sum g_l\o r_lx_{i_l}\in L(\Phi,R)$ and
$T(E(\Phi,R))\subseteq L(\Phi,R)$.

Let us prove the opposite inclusion. Clearly, $T(E(\Phi,R))$ contains the
nilpotent part $\{x_i\}$ of the Chevalley basis:
$x_i(t)=\exp(tx_i)=1+tx_i+o(t)$.  The remaining basis vectors $h_i$ lie in
$T(E(\Phi,R))$ also, because $h_i=x_i(1)\o x_{-i}+x_i-x_{-i}$ (see, e.g.,
[Bor70]). This completes the proof.

\medskip

In particular, Proposition 6 shows that any adjoint Chevalley group is
adjoint in the sense of Section 4.

\s 6.
Proof of the main theorems

{\bf The automorphisms of $E(\Phi,R)$.}
By Lemma 2 (vi), we have the natural injective homomorphism
$\Pi\:\Aut_\Z L(\Phi,R)\to\Aut E(\Phi,R)$. By Proposition 1 and Lemma 2
(i), each automorphism of $E(\Phi,R)$ is unipotent (for a suitably chosen
$m$) and, hence,  continuous (by Proposition 4 and Lemma~2 (i) and
(iv)) and, therefore, smooth by Lemma 4.  Then the map $\Pi$ is surjective
by Propositions 3 and 6. Thus, $\Aut E(\Phi,R)\iso\Aut_\Z L(\Phi,R)\iso
\Aut_\Z R\semitimes\Aut_R L(\Phi,R)$ (the latter isomorphism holds by
Lemma 2 (v)).

\noindent
{\bf The automorphisms of $G(\Phi,R)$} are the same as of $E(\Phi,R)$.
Indeed, each automorphism of $E(\Phi,R)$ is standard and, hence, can be
extended to an automorphism of $G(\Phi,R)$ by Lemma 2 (vi).  Thus, the
natural map $\Aut G(\Phi,R)\to \Aut E(\Phi,R)$ is surjective (and
well defined by Lemma 2 (iii)).  This map is also injective, because
of Lemma 2 (vi) and the following general fact.

\Lemma 5.
If $A$ is an automorphism invariant subgroup of a group $B$ and the
centraliser of $A$ in $B$ is trivial, then the natural map
$\rho\:\Aut B\to\Aut A$ is injective.

\Proof
For any $\phi\in\ker\rho$, $a\in A$, and $b\in B$, we have
$bab^{-1}=\phi(bab^{-1})=\phi(b)\phi(a)\phi(b^{-1})=\phi(b)a\phi(b^{-1})$.
Therefore, $b^{-1}\phi(b)$ centralises $A$.  Hence, $b=\phi(b)$ for any
$b\in B$. This completes the proof of the automorphism theorem.

\medskip

\noindent
{\bf The isomorphism theorem} is an easy corollary of the automorphism
theorem. Each isomorphism of Chevalley groups
$\sigma\:G(\Phi,R)\to G(\Phi,R')$ induces an automorphism $\phi_\sigma$ of
the group $G(\Phi,R\times R')$, because this group is isomorphic to
$G(\Phi,R)\times G(\Phi,R')$ and we can put
$\phi_\sigma(g,g')=(\sigma^{-1}(g'),\sigma(g))$. The standardness of
$\phi_\sigma$ implies that $\sigma$ is induced by an isomorphism of the
corresponding Lie rings. A similar argument applies to elementary
subgroups.

\s 7.
Automorphisms of Chevalley algebras

Recall that an \emph{inner automorphism} of a Chevalley algebra
$L(\Phi,R)$ is a conjugation $x\mapsto gxg^{-1}$ by an element $g$ of the
corresponding Chevalley group $G(\Phi,R)$. Clearly, the inner
automorphisms form a group isomorphic to $G(\Phi,R)$.

Let $\Delta=\{\delta_1,\dots,\delta_d\}$ be the symmetry group of
the Dynkin diagram of $\Phi$ (the number $d$ can be 1, 2, or 6, depending
on $\Phi$) and let $R=R_1\oplus\dots\oplus R_d$ be a (possibly trivial)
decomposition of the ring $R$ into a direct sum of ideals.  Suppose that
$f_i\in\Aut_{R_i}L(\Phi,R_i)$ is the automorphism induced by the symmetry
$\delta_i$ (see [VOn88]). The automorphism $f$ of the algebra
$L(\Phi,R)=L(\Phi,R_1)\oplus\dots\oplus L(\Phi,R_d)$ that sends
$x_1+\dots+x_d$ to $f_1(x_1)+\dots+f_d(x_d)$, where $x_i\in L(\Phi,R_i)$,
is called a \emph{diagram automorphism} of the algebra $L(\Phi,R)$.
Clearly, diagram automorphisms form a group isomorphic to the
subgroup
$$
D(\Phi,R)=\left\{\sum e_i\delta_i\;\mid\;
e_i\in R,\ e_i^2=e_i,\ e_ie_j=0 \hbox{ for } i\ne j,\ \sum e_i=1\right\}
$$
of the group of units of the group algebra $R\Delta$.

\Th 1.
Let $R$ be an associative commutative ring without additive torsion, with
unity and $1\over6$ and let $\Phi$ be a reduced irreducible root
system. Then any automorphism $f$ of the Lie $R$-algebra $L(\Phi,R)$ can
be expressed uniquely as a composition of diagram and inner automorphisms,
$\Aut_R L(\Phi,R)\iso D(\Phi,R)\semitimes G(\Phi,R)$.

\Proof
Let $n$ be the dimension of the Lie algebra $L(\Phi)$.
Consider the ideal $J$ in $\Z[x_{11},x_{12},\dots,x_{nn}]$ defining the
group $\Aut_\C L(\Phi)$. The ideal $J$ decomposes into a product
$J=J_1J_2\dots J_d$ of prime ideals $J_i$ corresponding to
irreducible (= connected) components
$h_iG(\Phi)$
of the group $\Aut_\C L(\Phi)$, where $h_i$ are integer matrices
of diagram automorphisms.
Take a matrix $A=(a_{pq})\in\Aut_R L(\Phi,R)$. Then
$f(a_{pq})=0$ for $f\in J$. Put $I_i=\{f(a_{pq})\;;\;f\in J_i\}\nin R$.
Then
\item{(i)}
$\prod I_i=\0$;
\item{(ii)}
$I_i+I_j=R$ for $i\ne j$ (otherwise we take the factor ring by a maximal
ideal $M\supseteq I_i+I_j$ and obtain a matrix $A_M$ belonging to the
intersection of two irreducible components of the group $\Aut_{R/M}
L(\Phi,R/M)$, but this intersection is empty, because $R/M$ is a field).

\enditem
These conditions (i) and (ii) imply that the ring $R$ is the direct sum
$R=\bigoplus R/I_i$ [Bou61, Ch.2 \S1, Proposition 5]. So, $A=\sum A_{I_i}$
and the entries of the matrix $A_{I_i}\in M_n(R/I_i)$ satisfy the
equations $f(a_{pq})=0$ for $f\in I_i$. Therefore,
$A_{I_i}=h_ig_i\in h_iG(\Phi, R/I_i)$ and $A=(\sum e_ih_i)(\sum g_i)$,
where $e_i$ is the unity of the ring $R/I_i$. This completes the proof.

Another approach to describing the automorphisms of Chevalley algebras
was suggested in [Pia02].

\REFERENCES

\baselineskip 11.7pt

\[Abe89]
Abe E.
Normal subgroups of Chevalley groups over commutative rings
{// Contemp. Math.}
1989. {V.83}. P.1--117.

\[Abe93]
Abe E.
Automorphisms of Chevalley groups over commutative rings
{// Algebra i Analiz.}
1993. {V.5}. no.2. P.74--90.

\[AHu88]
Abe E., Hurley J.
Centers of Chevalley groups over commutative rings
{// Comm. Algebra.}
1988. {V.16}. no.1. P.57--74.

\[ASu76]
Abe E., Suzuki K.
On normal subgroups of Chevalley groups over commutative rings
{// T\^ohoku Math. J.}
1976. {V.28}. no.1. P.185--198.

\[Bor70]
Borel A.
Properties and linear representations of Chevalley groups.
in Semin. Algebr. Groups related Finite Groups Princeton 1968/69, Lect.
Notes Math. 131, A1-A55 (1970).

\[Bou61]
Bourbaki N.
\'El\'ements de Math\'ematique. Alg\`ebre commutative. Chapitres 1 et 2.
Paris: Hermann. 1961.

\[Bun07]
Bunina E. I.
Automorphisms of elementary adjoint Chevalley groups of types
$A_l, D_l, E_l$ over local rings with~1/2
{// arXiv:math/0702046}. (2007). To appear in {Algebra i Logika}.

\[Che95]
Chen Yu
Automorphisms of simple Chevalley groups over $\Q$-algebras
{// T\^ohoku Math. J.}
1995. {V.47}. no.1. P.81--97.

\[Che96]
Chen Yu
Isomorphisms of adjoint Chevalley groups over integral domains
{// Trans. Amer. Math. Soc.}
1996. {V.348}. no.2. P.521--541.

\[Che00]
Chen Yu
Isomorphisms of Chevalley groups over algebras
{// J. Algebra.}
2000. {V.226}. no.2. P.719--741.

\[CKe99]
Costa D. L., Keller G. E.
On the normal subgroups of $G_2(A)$
{// Trans. Amer. Math. Soc.}
1999. {V.351}. no.12. P.5051--5088.

\[Gol97]
Golubchik I. Z.
Lie type groups over PI-rings
{// Fund. i Prikl. Mat.}
1997. {V.3}. no.2. P.399--424.

\[GMi83]
Golubchik I. Z., Mikhalev A. V.
Isomorphisms of unitary groups over associative rings.
{// Zap. Nauchn. Sem. Leningr. Otd. Mat. Inst. Akad. Nauk SSSR.}
1983. {V.132}. P.97--109.

\[HaV03]
Hazrat R., Vavilov N.
$K_1$ of Chevalley groups are nilpotent
{//J. Pure Appl. Algebra.} 2003. {V.179}. P.99--116.

\[HO'M89]
Hahn A. J, O'Meara O. T.
The classical groups and K-theory.
Springer. Berlin et al. 1989.

\[Lev83]
Levchuk V.M.
Connections between a unitriangular group and certain rings.
Chap. 2: Groups of automorphisms
{// Sibirsk. Matem. Zh.} 1983. {V.24}. no.4. P.543--557

\[Pet82]
Petechuk V. M.
Automorphisms of matrix groups over commutative rings
{// Mat. Sb.}
1982. {V.117}. no.4. P.534--547.

\[Pia02]
Pianzola A.
Automorphisms of toroidal Lie algebras their central quotients
{// J. Algebra and Appl.}
2002. {V.1}. no.1. P.113--121.

\[Vas86]
Vaserstein L. N.
On normal subgroups of Chevalley groups over commutative rings
{// T\^ohoku Math. J.}
1986. {V.36}. no.5. P.219--230.

\[VGN06]
Vavilov N. A., Gavrilovich M. R., Nikolenko S. I.
The structure of Chevalley groups: a Proof from the book
{// Zap. Nauchn. Sem. S-Peterburg. Otd. Mat. Inst. Steklov.}
2006. {V.330}. no.13. P.36--76.

\[VPl96]
Vavilov N., Plotkin E.
Chevalley groups over commutative rings.
I: Elementary calculations.
{// Acta Appl. Math.}
1996. {V.45}. no.1. P.73--113.

\[VOn88]
Vinberg E.B., Onishchik A.L.
Seminar on Lie groups and algebraic groups.
Moscow: Nauka. 1988.

\[Wat80]
Waterhouse W. C.
Automorphisms of GL$_n(R)$
{// Proc. Amer. Math. Soc.}
1980. {V.79}. no.3. P.347--351.

\[Zal83]
Zalesskii A. E.
Linear groups
// Itogi nauki i tekhn. VINITI. Algebra. Topologiya. Geometriya.
1983. Moscow. V.21. P.135--182.

\[Zel85]
Zelmanov E. I.
Isomorphisms of linear groups over an associative ring
{// Sibirsk. Matem. Zh.}
1985. {V.26}. no.4. P.49--67.

\end